\documentclass[12pt]{article}
\usepackage{theorem,amsfonts,amssymb}
\newtheorem{theorem}{Theorem}[section]
\newtheorem{proposition}{Proposition}[section]
\newtheorem{lemma}{Lemma}[section]
\newtheorem{corollary}{Corollary}[section]
\theorembodyfont{\upshape}
\newtheorem{definition}{Definition}
\newtheorem{remark}{Remark}[section]
\newtheorem{example}{Example}[section]
\newtheorem{proof}{Proof}

\newcommand{\bt}{\begin{theorem}}
\newcommand{\et}{\end{theorem}}
\newcommand{\bl}{\begin{lemma}}
\newcommand{\el}{\end{lemma}}
\newcommand{\bp}{\begin{proposition}}
\newcommand{\ep}{\end{proposition}}
\newcommand{\bex}{\begin{example}}
\newcommand{\eex}{\end{example}}
\newcommand{\bc}{\begin{corollary}}
\newcommand{\ec}{\end{corollary}}
\newcommand{\bo}{\begin{proof}}
\newcommand{\eo}{\end{proof}}
\newcommand{\bd}{\begin{definition}}
\newcommand{\ed}{\end{definition}}
\newcommand{\br}{\begin{remark}}
\newcommand{\er}{\end{remark}}
\newcommand{\be}{\begin{enumerate}}
\newcommand{\ee}{\end{enumerate}}

\begin{document}

\title{Operator decomposable measures and stochastic difference equation}
\author{C. R. E. Raja}
\date{}
\maketitle

\let\ol=\overline
\let\epsi=\epsilon
\let\vepsi=\varepsilon
\let\lam=\lambda
\let\Lam=\Lambda
\let\ap=\alpha
\let\vp=\varphi
\let\ra=\rightarrow
\let\Ra=\Rightarrow
\let \Llra=\Longleftrightarrow
\let\Lla=\Longleftarrow
\let\lra=\longrightarrow
\let\Lra=\Longrightarrow
\let\ba=\beta
\let\ov=\overline
\let\ga=\gamma
\let\Ba=\Delta
\let\Ga=\Gamma
\let\Da=\Delta
\let\Oa=\Omega
\let\Lam=\Lambda
\let\ups=\upsilon

\newcommand{\Ad}{{\rm Ad}}
\newcommand{\cK}{{\cal K}}
\newcommand{\pc}{T_{si}}
\newcommand{\Aut}{{\rm Aut}}
\newcommand{\cR}{{\cal R}}
\newcommand{\Var}{{\rm Var}}
\newcommand{\cG}{{\cal G}}
\newcommand{\cH}{{\cal H}}
\newcommand{\G}{{\mathbb G}}
\newcommand{\Z}{{\mathbb Z}}
\newcommand{\Q}{{\mathbb Q}}
\newcommand{\cN}{{\cal N}}
\newcommand{\cE}{{\cal E}}
\newcommand{\N}{{\mathbb N}}
\newcommand{\R}{{\mathbb R}}
\newcommand{\C}{{\mathbb C}}
\newcommand{\T}{{\mathbb T}}
\newcommand{\mP}{{\mathbb P}}



\begin{abstract}
We consider the following convolution equation or equivalently
stochastic difference equation
$$\lam _k = \mu _k*\phi (\lam _{k-1}),~~~~~~ k \in \Z \eqno (1) $$
for a given bi-sequence $(\mu _k)$ of probability measures on $\R
^d$ and a linear map $\phi $ on $\R ^d$.  We study the solutions of
equation (1) by realizing the process $(\mu _k)$ as a measure on
$(\R ^d)^\Z$ and rewriting the stochastic difference equation as
$\lam = \mu
*\tau (\lam )$-any such measure $\lam$ on $(\R ^d)^\Z$ is known
as $\tau$-decomposable measure with co-factor $\mu$-where $\tau$ is
a suitable weighted shift operator on $(\R ^d)^\Z$.  This enables
one to study the solutions of (1) in the settings of
$\tau$-decomposable measures.  A solution $(\lam _k)$ of (1) will be
called a fundamental solution if any solution of (1) can be written
as $\lam _k*\phi ^k(\rho )$ for some probability measure $\rho $ on
$\R ^d$. Motivated by the splitting/factorization theorems for
operator decomposable measures, we address the question of existence
of fundamental solutions when a solution exists and answer
affirmatively via a one-one correspondence between fundamental
solutions of (1) and strongly $\tau$-decomposable measures on $(\R
^d)^\Z$ with co-factor $\mu$.  We also prove that fundamental
solutions are extremal solutions and vice versa.  We provide a
necessary and sufficient condition in terms of a logarithmic moment
condition for the existence of a (fundamental) solution when the
noise process is stationary and when the noise process has
independent $\ell _p$-paths.
\end{abstract}

\begin{section}{Introduction}

Stochastic/random difference equations arise in different contexts
and studied by many since \cite{Ke}, \cite{Ts} and \cite{Yo}.  We
consider the following stochastic difference equation

$$\eta _k = \xi _k + \phi ( \eta _{k-1} ), ~~~~ k \in \Z \eqno {\rm(TY)} $$
where $\eta _k$ and $\xi _k$ are $\R ^d$-valued random variables and
$\phi$ is a linear map on $\R^d$ (with the assumption that for each
$k\in \Z$, $\{ \xi _k, \eta _{k-i} \mid i \leq k \}$ are
independent).  The random variables $(\xi _k)$ are given and are
called the noise process of the equation (TY). The random variables
$(\eta _k)$ satisfying equation (TY) for given $(\xi _k)$ and $\phi$
will be called a solution of equation (TY).  We wish to study the
distributional properties of solutions of equation (TY). This pushes
us to focus on the corresponding convolution equation

$$\lam _k = \mu _k * \phi (\lam _{k-1}), ~~~~~~ k \in \Z \eqno (1) $$
where $\lam _k$ and $\mu _k$ are the laws of $\eta _k$ and $\xi_k$
($k \in \Z$) respectively.  It may be noted that Lemma 4.3(ii) of
\cite{HiY} asserts that for a solution $(\lam _k)$ of equation (1)
there exists a solution $(\eta _k)$ of equation (TY) whose marginal
laws are $(\lam _k)$.

B. Tsirelson \cite{Ts} considered the following stochastic
difference equation on the real line
$$\eta _k = \xi _k +{\rm frac} (\eta _{k-1} ) ~~~~ k \in -\N
\eqno (2)$$ with a given stationary Gaussian noise process $(\xi_k)$
to obtain his celebrated example of the stochastic differential
equation
$$dX_t = dB_t+ b^+(t, X)dt, ~~ X(0)=0 $$
that has unique but not strong solution where  $(B_t)$ is the
one-dimensional Brownian motion (see \cite{Ts} for more details). M.
Yor \cite{Yo} formulated Tsirelson's equation (2) in the form of
equation (TY) on the one-dimensional torus $\R /\Z$ when $\phi$ is
the identity automorphism and $(\xi _k)$ is a general noise process.
In particular, \cite{Yo} identified solutions of equation (TY) with
measures on a quotient $(\R/\Z)/M$ where $M$ is a closed subgroup of
$\R/\Z$.  When $\phi$ is the identity automorphism, equation $(1)$
was considered on general compact groups \cite{AUY} and when the
noise law $(\xi _k)$ is stationary, equation $(1)$ was considered on
abelian groups \cite{Ta}.  When the noise law is stationary,
\cite{Ra} considered equation (1) on general locally compact groups
and provided a necessary and sufficient condition for existence of
solutions, under further distal condition \cite{Ra} completely
classified solutions of equation (1).

Here, we study equation (1) by realizing the stochastic processes
$(\eta _k)$ and $(\xi _k)$ in the space $(\R ^d)^\Z$.  Recall that
$(\R^d)^\Z$ is the space of all bi-sequences $(v_n)_{n \in \Z}$ of
vectors in $\R ^d$ and $(\R ^d)^\Z$ is a complete separable metric
vector space (of infinite dimension).  The metric $d$ on $(\R^
d)^\Z$ is given by
$$d((v_n), (w_n)) = \sup _{n \in\Z} \min \{ ||v_n-w_n||, 1/|n| \}$$
and with this metric a sequence in $(\R ^d)^\Z$ converges if and
only if the sequence of each co-ordinates converges.

For any $k_1, \cdots , k_n$ in $\Z$, let $P_{k_1, \cdots , k_n}
\colon (\R ^d)^\Z \to (\R^d)^n$ be the projection defined by
$P_{k_1, \cdots k_n}(x_k) = (x_{k_1}, \cdots , x_{k_n})$.  Then
$P_{k_1, \cdots , k_n}$ is a continuous linear map.

We now look at the realization of $(\lam _k)$ and $(\mu _k)$ (that
satisfy equation (1)) in $(\R ^d)^\Z$.  Using Kolmogorov consistency
theorem we can find (unique) probability measures $\lam$ and $\mu$
on $(\R ^d)^\Z$ such that $$P_{k_1, \cdots , k_n}(\lam ) = \lam
_{k_1}\times \cdots \times \lam _{k_n} ~~{\rm and}~~ P_{k_1, \cdots
, k_n}(\mu ) = \mu _{k_1}\times \cdots \times \mu _{k_n}$$ for any
set of finite integers $k_1, \cdots , k_n$.  In particular,
$P_k(\lam ) = \lam _k$ and $P_k(\mu ) = \mu _k$ for all $k \in \Z$.
In this situation, we sometime denote $\mu$ also by $(\mu _k)$.

In order to rephrase equation (1) in $(\R ^d)^\Z$, we need the
following weighted shift operator.   Define $\tau \colon (\R ^d)^\Z
\to (\R ^d)^\Z $ by
$$\tau (v_i) = (\phi (v_{i-1})) ,~~~~(v_i)\in (\R ^d)^\Z \eqno
{\rm (L)}$$ which is the composition of the right shift operator and
the diagonal action of $\phi$ . It can easily be seen that $\tau$ is
a continuous linear operator on the complete separable metric vector
space $(\R ^d)^\Z$ and $\tau$ is invertible if and only if $\phi$ is
invertible.  With the help of the weighted shift operator $\tau$ and
the uniqueness part of Kolmogorov consistency theorem, equation (1)
can be rewritten as the following convolution equation on $(\R
^d)^\Z$

$$\lam = \mu *\tau (\lam). \eqno (3)$$

Equation (3) motivates us to look at the so called operator
decomposable measures on vector spaces studied by Siebert et al
(\cite{Si1}, \cite{Si2} and the references cited therein).

We now recall operator decomposable measures on complete separable
metric vector spaces. Let $V$ be a complete separable metric vector
space and $M^1(V)$ be the space of Borel probability measures on
$V$.  Given a continuous linear operator $T$ on $V$, $\rho \in M^1(V
)$ is called $T$-decomposable if there is a $\nu \in M^1(V)$ such
that $$\rho = \nu
*T(\rho )$$ and in this situation we also say that $\rho$ is
$T$-decomposable with co-factor $\nu$.

In the language of decomposable measures, a solution $(\lam _k)$ of
equation (1) gives a $\tau$-decomposable measure $\lam$ on $(\R
^d)^\Z$ with co-factor $\mu$ such that $P_k(\lam ) = \lam _k$ and
$P_k(\mu ) = \mu _k$ for all $k \in \Z$.  We use this approach to
study the distributional properties of solutions of the stochastic
difference equation (TY).   For instance, this approach is useful to
understand when solutions of equation (1) can have atoms: it may be
recalled that if $\mu _k$ has density, then $\lam _k$ also has
density that is even smoother.

\bp\label{atom} Let $(\mu _k)$ be a bi-sequence in $M^1(\R ^d)$ and
$\phi$ be a linear map on $\R ^d$.  Assume that $(\lam _k)$ is a
solution of equation (1). If $\inf _k \prod _{i=-k}^k \lam _i(\{ x_i
\}) \not =0$ for some sequence $(x_k)_{k\in \Z}$ in $\R ^d$, then
each $\mu _k$ is a dirac measure for all $k$. \ep

The above condition that  $\inf _k \prod _{i=-k}^k \lam _i(\{ x_i
\}) \not =0$ may be viewed as  $\lam _k(\{ x_k \}) \to 1 $ faster as
$|k|\to \infty$.  Thus, conclusion of Proposition \ref{atom} may be
read as $\lam _k(\{ x_k \}) $ does not converge to $1$ faster as
$|k|\to \infty$ for any sequence $(x_k)$ unless every $\mu _k$
degenerates.

Among $T$-decomposable measures, strongly $T$-decomposable measures
(that is, $T$-decomposable measure $\rho$ with $T^n(\rho ) \to
\delta _0$) and $T$-invariant measures (that is, $T(\rho ) = \rho
*\delta _x$ for some $x\in V$) are particular cases.  Significance
of these two particular classes of $T$-decomposable measures stems
from the factorization theorem of Siebert \cite{Si2}: recall that
factorization theorem of \cite{Si2} proves that any symmetric
$T$-decomposable measure on a separable Banach space is a product of
a $T$-invariant measure and a strongly $T$-decomposable measure and
a similar result was proved for $T$-decomposable measures verifying
certain nondegeneracy condition in \cite{Si1}.

This motivates us to look for similar interesting classes among
solutions of equation (1).  At this time we note that $(\lam _k
*\phi ^k(\rho ))$ is also a solution of (1) if so is $(\lam _k)$.
In view of these reasons a solution $(\lam _k)$ of (1) will be
called a fundamental solution if to each solution $(\nu _k)$ of
equation (1), there is a $\nu \in M^1(\R ^d)$ such that $\nu _k=
\lam _k *\phi ^k(\nu )$ for all $k \in \Z$.  Before we proceed to
establish the correspondence between fundamental solutions and
strongly $\tau$-decomposable measures, we would look at the other
component of the factorization theorem of \cite{Si2}, the so-called
strictly $\tau$-invariant measure $\lam \in M^1((\R ^d)^\Z)$ (that
is, $\tau (\lam ) = \lam$). If $\lam \in M^1((\R ^d)^\Z)$ is
strictly $\tau$-invariant, then it is easy to see that for any $k
\in \Z$, $P_k(\lam ) = \phi ^k(P_0(\lam ))$.  Thus, forced by the
factorization theorem of \cite{Si2}, we pose the question: is there
any fundamental solution of equation (1) if equation (1) has a
solution.  We provide an affirmative answer via a one-one
correspondence between fundamental solutions and strongly
$\tau$-decomposable measures.

\bt\label{sdfs}  Let $(\mu _k)_{k\in \Z}$ be a bisequence in $M^1(\R
^d)$ and $\phi$ be an invertible linear map of $\R ^d$.  Let $\tau$
be the linear operator on $(\R ^d)^\Z$ defined by the equation (L).

\be
\item [(i)] If $\lam \in M^1((\R ^d)^\Z)$ is a strongly $\tau$-decomposable
measure with co-factor $\mu = (\mu _k)\in M^1((\R ^d)^\Z)$, then
$\lam _k = P_k(\lam )$ is a fundamental solution of equation (1).

\item [(ii)] If equation (1) has a solution, then there is a strongly
$\tau$-decomposable measure on $(\R ^d)^\Z$ with co-factor $\mu
*\delta _v$ for some $v\in (\R ^d)^\Z$ where $\mu = (\mu _k)\in
M^1((\R ^d)^\Z)$. In particular, equation (1) has a fundamental
solution. \ee \et

For $\mu \in M^1(V)$, it is easy to verify that if there is a
strongly $T$-decomposable measure with co-factor $\mu$, then it is
unique (see 1.2 and 1.3 of \cite{Si1}). Since fundamental solution
is defined so as to have similarities with strongly
$\tau$-decomposable measures, we ask the question that could there
be two fundamental solutions: before the answer is stated, it is
worth to note that $(\lam _k*\phi ^k(\delta _v))$ is also a
fundamental solution if so is $(\lam _k)$.

\bp\label{ufs} If $(\lam _k)$ and $(\rho _k)$ are fundamental
solutions of (1), then there is a $v\in \R ^d$ such that $\rho _k=
\lam _k *\phi ^k(\delta _v)$ for all $k \in \Z$.\ep

It may be easily seen that the set of all solutions of equation (1)
is a convex set and any extreme point of this convex set is known as
extremal solution of (1).  We now look at the relation between
extremal solutions and fundamental solutions.

\bt\label{fes} Let $(\mu _k)$ and $\phi$ be as in Theorem
\ref{sdfs}. Then $(\lam _k)$ is a fundamental solution of equation
(1) if and only if $(\lam _k)$ is a extremal solution.  In
particular, the set of extremal solutions is either empty or can be
identified with $\R ^d$.\et

Thus, it follows from Theorem \ref{fes} and Proposition \ref{ufs}
that having one extremal solution is sufficient to get all
(extremal) solutions.

Siebert \cite{Si1} provides a necessary and sufficient logarithmic
moment condition for a measure on Banach space to be a co-factor of
strongly operator decomposable measure by generalizing a result of
Zakusilo \cite{Za}.  Motivated by this result we provide a similar
condition on the noise process for the existence of a (fundamental)
solution of equation (1) when the noise is stationary (Theorem
\ref{sn}) and when the noise has independent $\ell _p$-paths
(Theorem \ref{bp}).

\end{section}

\begin{section}{Preliminaries}

We consider a complete separable metric vector space $V$.  Let
$M^1(V)$ denote the space of all Borel probability measures on $V$
with weak topology with respect to all continuous bounded functions
on $V$.

One of the useful fact about the weak topology is the criterion for
compact sets in $M^1(V)$ given by the following:

\noindent {\bf Prokhorov's Theorem (cf. Chapter 2, Theorem 6.7 of
\cite{Pa}):} {\it A subset $\Ga$ of $M^1(V)$ has compact closure or
equivalently relatively compact if and only if for each $\epsilon
>0$ there is a compact set $K_\epsilon \subset V$ such that $\nu
(K_\epsilon )
>1-\epsilon$ for any $\nu \in \Ga$.}

For a measure $\nu \in M^1(V)$, we define $\check \nu \in M^1(V)$ by
$\check \nu (E) = \nu (-E)$ for any Borel subset $E$ of $V$.  A
measure $\nu \in M^1(V)$ is called symmetric if $\check \nu = \nu$.

For any $x\in V$, $\delta _x \in M^1(V)$ denotes the dirac measure
supported at $\{ x\}$.  For any $\nu _1 , \nu _2 \in M^1(V)$, we
denote the convolution product of $\nu _1$ and $\nu _2$ by $\nu
_1*\nu _2$ and is defined by
$$\nu _1* \nu _2 (E) = \int \nu _2 (x+E) d\nu _1(x)$$ for any Borel
subset $E$ of $V$.  The convolution product $\delta _x*\nu $ is
known as shift of $\nu \in M^1(V)$ by $x\in V$.

We say that a measure $\nu _1\in M^1(V)$ is a factor of a measure
$\nu \in M^1(V)$ if there exists a $\nu _2 \in M^1(V)$ such that
$\nu = \nu _1*\nu _2$.

For a continuous linear operator $T$ on $V$ and $\nu \in M^1(V)$, we
define $T(\nu ) \in M^1(V)$ by $T(\nu )(E) = \nu (T^{-1}(E))$ for
any Borel subset $E$ of $V$.  Then it can easily be seen that $T(\nu
*\nu ') = T(\nu )*T(\nu ')$ for any $\nu , \nu ' \in M^1(V)$.

We will be studying solutions $(\lam _k)$ of the convolution
equation
$$ \lam _k = \mu _k *\phi (\lam _{k-1}), ~~~~~ k\in \Z \eqno(1) $$
when $V= \R ^d$ and $\phi $ is a linear transformation on $\R ^d$.

We often make use of the following handy application of Prokhorov's
theorem to convolution equations.

\noindent{\bf Theorem (cf. Chapter 3 Section 2 of \cite{Pa}):} {\it
Let $(\rho _n), (\sigma _n)$ and $(\nu _n)$ be sequences in $M^1(V)$
such that $\rho _n = \sigma _n *\nu _n $ for all $n \geq 1$.

\be \item[(a)] If $(\rho _n)$ is relatively compact, then there is a
sequence $(x_n)$ in $V$ such that $(\sigma _n *\delta _{x_n})$ is
relatively compact.

\item[(b)] If two of the three sequences $(\rho _n), (\sigma _n)$ and
$(\nu _n)$ are relatively compact, then so is the third one. \ee}

For symmetric measures we have the following improvement of the
above, proof which is essentially based on methods in 3.3 of
\cite{Si2}.

\bl\label{sym} Let $(\rho _n)$, $(\nu _n)$ and $(\sigma _n)$ be
sequences in $M^1(\R ^d)$.  Suppose $\rho _n = \sigma _n
*\nu _n$ for all $n$ and $(\rho _n)$ is relatively compact.  If each
$\sigma _n$ is symmetric and $\sigma _n$ is a factor of $\sigma_m$
for all $n\geq m$ or for all $m \geq n$, then $(\sigma _n)$
converges in $M^1(\R ^d)$.\el

\bo By Corollary 2.5.3 of \cite{Li}, we get that $(\sigma _n )$ is
relatively compact in $M^1(\R ^d)$.  If $\sigma _n = \sigma
_{m}*\varrho _{n,m}$ for all $n ,m$ with $n>m$ for some $\varrho
_{n,m}\in M^1(\R ^d)$, let $\sigma$ and $\sigma '$ be two limit
points of $(\sigma _n)$. Let $(l_m)$ be a subsequence such that
$\sigma _{l_m}\to \sigma$.  Fix $n \geq 1$ and consider $\sigma
_{l_m} = \sigma_{n}
* \varrho _{l_m,n}$ for all $l_m\geq n$.  This implies that
$(\varrho _{l_m,n})$ is relatively compact.  Then there is a
$\varrho _n'\in M^1(\R ^d)$ such that $\sigma = \sigma _n*\varrho
_n'$ for all $n$. This implies that $\sigma = \sigma '*\varrho $ for
some $\varrho \in M^1(\R ^d)$. Similarly, we get that $\sigma '=
\sigma *\varrho '$ for some $\varrho '\in M^1(\R ^d)$.  This implies
that $\sigma ' = \sigma '*\varrho '*\varrho $, hence $\varrho
'*\varrho = \delta _0$. Thus, $\varrho = \delta _x$, hence $\sigma =
\sigma ' *\delta _x$ for some $x\in V$.  Since each $\sigma _n$ is
symmetric, the limit points $\sigma $ and $\sigma '$ of $(\sigma
_n)$ are also symmetric.  But $\sigma = \sigma ' *\delta _x$, hence
$\sigma = \sigma '$.  Thus, $(\sigma _n )$ converges. Convergence of
$(\sigma _n)$ can be proved in a similar way if $\sigma _n$ is a
factor of $\sigma_m$ for all $n \geq m$. \eo

We also need a reformulation of Theorem 3.1 of \cite{Cs} and this
reformulation is proved by considering semidirect products. Given a
continuous invertible linear operator $T$ on a complete separable
metric vector space $V$, the semidirect product of $\Z$ and $V$ with
respect to $T$ is denoted by $\Z \ltimes _T V$ whose underlying
space is $\Z \times V$ and the group multiplication is given by:
$(n,v)(m,w)= (n+m, v+T^n(w))$ for $n,m \in \N$ and $v,w \in V$. Then
$\Z \ltimes _T V$ is a complete separable metric group.  For any
measure $\mu \in M^1(V)$ and $n \in \Z$, we define $n \otimes \mu$
by $(n\otimes \mu)(A\times B) = \delta _n(A) \mu (B)$ for any subset
$A$ of $\Z$ and any Borel subset $B$ of $V$.  Then $n\otimes \mu$ is
a probability measure on $\Z \ltimes _T V$.  Since $v\mapsto (0,v)$
is an isomorphism of $V$ onto its image, $V$ will be realized as a
subgroup of $G$, hence any measure on $V$ is also realized as a
measure on $G$.  Thus, $0\otimes \mu$ will be simply denoted by
$\mu$ for any measure $\mu \in M^1(V)$.  We now present a
reformulation of Theorem 3.1 of \cite{Cs} suitable for our study of
equation (1).

\bl\label{sft} Let $(\mu _k)$ be a bi-sequence in $M^1(\R ^d)$ and
$\phi$ be an invertible linear map on $\R ^d$.  Suppose equation (1)
has a solution. Then there is a sequence $(w_n)$ in $\R ^d$ such
that the sequence $(\prod _{i=0}^{n-1} \phi ^i(\mu _{k-i}) * \phi
^{n}(\delta _{w_{n-k}}) )$ converges for any $k <0$. \el

\bo  If there is a solution $(\nu _k)$ of equation (1), then $\nu
_{-1} = \mu _{-1}*\phi (\nu _{-2}) = \cdots = \prod _{i=0}^{n-1}\phi
^i(\mu _{-i-1}) *\phi ^n(\nu _{-1-n})$.  Then there is a $(a_n)$ in
$\R ^d$ such that $(\prod _{i=0}^{n-1}\phi ^i(\mu _{-i-1})
*\delta _{a_n})$ is relatively compact.

Let $G= \Z \ltimes _\phi \R ^d$ and for $n \geq 1$, let $\lam _n = 1
\otimes \mu _{-n} \in M^1(G)$.  Then $$\lam _1 *\cdots * \lam _n
*\delta _{(-n, \phi ^{-n}(a_n))} = \mu _{-1} * \phi (\mu _{-2}) *\cdots *\phi
^{n-1}(\mu _{-n})*\delta _{a_n}$$ for all $n$.   By Theorem 3.1 of
\cite{Cs}, there exist $w_n$ in $\R ^d$, $k_n \in \Z$ and $\nu _k
\in M^1(G)$ such that $\lam _{k+1} *\cdots *\lam _n *\delta
_{(k_{n},w_{n})} \to \nu _k$ as $n \to \infty$ for $k\geq 0$. This
implies for each $k <0$ that there are $\varrho _k \in M^1(\R ^d)$
such that $\mu _k*\phi (\mu _{k-1})
*\cdots *\phi ^{k+n-1}(\mu _{-n+1})*\phi ^{k+n}(\delta _{w_{n}}) \to
\varrho _k$ or equivalently $\prod _{i=0}^{n-1} \phi ^i(\mu _{k-i})
* \phi ^{n}(\delta _{w_{n-k}}) \to \varrho _k$ as $n \to \infty$.
\eo

As our approach involves operator decomposable measures, we recall
the following:

\bd For a continuous linear operator $T$ on $V$, a measure $\rho \in
M^1(V)$ is called $T$-decomposable if there is a $\nu \in M^1(V)$
such that $\rho = \nu *T(\rho )$ and in this case $\nu$ is called
co-factor. A $T$-decomposable $\rho $ is called strongly
$T$-decomposable if $T^n (\rho ) \to \delta _0$.  \ed

We now recall the other component of splitting/factorization
Theorems of \cite{Si1} and \cite{Si2}.

\bd For a continuous linear operator $T$ on $V$, a measure $\rho \in
M^1(V)$ is called $T$-invariant (resp. strictly $T$-invariant) if
there is a $x\in V$ such that $\rho = T(\rho )*\delta _x$ (resp.
$T(\rho ) = \rho$).\ed

As explained in the introduction, by realizing the $\R^d$-valued
stochastic processes in the space $(\R^d) ^\Z$ of bi-sequences of
vectors, we study the solutions to equation (1). It is easy see that
$(\lam _k)$ is a solution to equation (1) if and only if the measure
$\lam $ is $\tau$-decomposable with co-factor $\mu$ where $\lam$ and
$\mu$ are given by Kolmogrov consistency theorem with marginal laws
given by $\lam _k$ and $\mu _k$ respectively and $\tau (x_k) = (\phi
(x_{k-1}))$.  We first obtain the proof on the atoms of solution of
equation (1).

\bo $\!\!\!\!\!${\bf ~of Proposition \ref{atom}~~~} By Kolmogrov
consistency theorem and by its uniqueness part there are $\lam$ and
$\mu$ in $M^1((\R ^d)^\Z)$ such that $P_k(\lam ) = \lam _k$,
$P_k(\mu )=\mu _k$ on $\R^d$ and $\lam = \mu *\tau (\lam )$ where
$\tau (v_i) = (\phi (v_{i-1}))$.  It follows from Lemma 3.5 of
\cite{Si1} that $\lam (\{ (x_k) \}) =0$ for any $(x_k)$ or $\mu$ is
a dirac measure. Now the result follows from
$$\lam (\{ (x_k) \}) = \lim _{k\to \infty} \prod _{i=-k}^k \lam
_i(\{ x_i \}) = \inf _k \prod _{i=-k}^k \lam _i(\{ x_i \}),
~~~~(x_k) \in (\R ^d)^\Z .$$ \eo

As observed in the introduction, for any solution $(\lam _k)$,
$(\lam _k*\phi ^k(\rho))$ is also a solution for $\rho \in M^1(\R
^d)$.  This motivates us to make the following definition:

\bd A solution $(\lam _k)$ of equation (1) is called a fundamental
solution of (1) if to each solution $(\nu _k)$ there is a $\nu \in
M^1(\R ^d)$ such that $\nu _k= \lam _k *\phi ^k(\nu )$ for all $k
\in \Z$.\ed

As we also consider infinitely divisible measures and process with
independent $\ell _p$-paths, we make the following formal
definitions:

\bd A measure $\rho \in M^1(\R ^d)$ is said to be infinitely
divisible if to each $n \in \N$, there is a $\rho _n \in M^1(\R ^d)$
such that $\rho _n^n = \rho$.\ed

\bd A $\R ^d$-valued stochastic process $(Y_k)$ (resp. a sequence of
probability measures $(\rho _k)$ on $\R ^d$) is said to have
independent $\ell _p$-paths for some $p \in [1, \infty]$ if there is
a $\R ^d$-valued stochastic process $(Y_k')$ such that $Y_k$ and
$Y_k'$ have the same law (resp. law of $Y_k'$ is $\rho _k$),
$(Y_k')$ is an independent bi-sequence and $(Y_k'(\omega ))$ is in
$\ell _p$ a.s.\ed

\end{section}

\begin{section}{Fundamental solution}

We now prove the results on fundamental solution.  We first provide
a useful sufficient condition for the existence of strongly
$\tau$-decomposable measures (or fundamental solutions) using
methods of tail idempotents of \cite{Cs}.

\bp\label{suff} Let $(\mu _k)$ be a bi-sequence of probability
measures on $\R ^d$ and $\phi$ is a linear map on $\R ^d$.  Define
the linear operator $\tau$ on $(\R ^d)^\Z$ by equation (L).  Suppose
$\mu _{k,n}= \mu _k* \phi (\mu _{k-1})*\cdots
*\phi ^n(\mu _{k-n})\to \rho _k \in M^1(\R ^d)$ as $n \to \infty$
for all $k$.  Then $(\rho _k)\in M^1((\R ^d)^\Z)$ is strongly
$\tau$-decomposable with co-factor $\mu = (\mu _k) \in M^1( (\R
^d)^\Z)$. \ep

\bo  By Kolmogorov consistency theorem there exists $\mu , \rho \in
M^1((\R ^d)^\Z)$ such that $P_k(\rho ) = \rho _k$ and $P_k(\mu )=
\mu _k$ for all $k \in \Z$.

For any $k\in \Z$ and $0\leq j \leq n$, we have $\mu _{k,n}= \mu
_{k,j-1}*\phi ^{j}(\mu _{k-j,n-j})$. Letting $n \to \infty$, we get
that $\rho _k = \mu _{k,j-1}*\phi ^{j}(\rho _{k-j})$ for any $j\geq
1$ and $k\in \Z$.  For $j=1$, $\rho _k = \mu _k *\phi (\rho _{k-1})$
for all $k \in \Z$.  By the uniqueness part in Kolmogorov
consistency Theorem, we get that $\rho = \mu *\tau (\rho )$ as
$P_k\tau = \phi P_{k-1}$. Thus, $\rho$ is $\tau$-decomposable with
co-factor $\mu$.

Using  $\rho _k = \mu _{k,j-1}*\phi ^{j}(\rho _{k-j})$ for any
$j\geq 1$ and $k\in \Z$, we get that $(\phi ^j(\rho _{k-j}))$ is
relatively compact. If $\rho _{k,\infty}$ is a limit point of $(\phi
^j(\rho _{k-j}))$, we have $\rho _k = \rho _k *\rho _{k, \infty}$.
This shows that $\rho _{k, \infty} = \delta _0$. Thus, $\phi ^j(\rho
_{k-j} )\to \delta _0$ as $j \to \infty$, hence $\tau ^j(\rho ) \to
\delta _0$ as $j \to \infty$.  Thus, $\rho$ is strongly
$\tau$-decomposable with co-factor $\mu$.  \eo

The following lemma is useful for two reasons: one reason is that it
gives a solution of equation (1) in case each $\mu _k$ degenerates
and other reason is that it explains why we have only strictly
$\tau$-invariant measures in the factorization of solutions to
equation (1), that is in the concept of fundamental solutions.

\bl\label{sdc}  Let $(x_k)$ be a given bi-sequence in $\R ^d$ and
$\phi$ be an invertible linear transformation on $\R ^d$.  Then
there is a $(y_k)$ such that $y_k = x_k+\phi (y_{k-1})$ for all $k
\in \Z$. In other words, for each $x \in (\R ^d)^\Z$ there is a
$y\in (\R ^d)^\Z$ such that $y= x+\tau (y)$.  \el

\bo Fix $y_0 \in \R ^d$.  Define $y_k = \phi ^k(y_0)+\sum _{i=1}^k
\phi ^{k-i}(x_i)$ for $k>0$.  Then $y_{k+1}= \phi ^{k+1}(y_0)+\sum
_{i=1}^{k+1} \phi ^{k+1-i}(x_i) = \phi (\phi ^k(y_0)+ \sum _{i=1}^k
\phi ^{k-i}(x_i)) +x_{k+1}= x_{k+1}+\phi (y_k)$ for all $k\geq 0$.

Define $y_k = \phi ^k(y_0)- \sum _{i=0}^{-k-1} \phi ^{k+i}(x_{-i})$
for $k <0$. Then $y_{k-1} = \phi ^{k-1}(y_0) -\sum _{i=0} ^{-k} \phi
^{k-1+i}(x_{-i}) = \phi ^{-1}(\phi ^k(y_0)-\sum _{i=0} ^{-k-1} \phi
^{k+i}(x_{-i}))-\phi ^{-1}(x_k)= \phi ^{-1}(y_k-x_k)$, hence $y_k =
x_k+\phi (y_{k-1})$ for all $k\leq 0$. \eo

The next result enables us to work with shifted noise process.

\bl\label{seqn} Let $(\mu _k)$ be a bi-sequence in $M^1(\R ^d)$ and
$\phi $ be an invertible linear map on $\R ^d$.  If $(\rho _k)$ is a
solution of equation of (1) with noise $(\mu _k)$, then to each
bi-sequence $(v_k)$ in $\R ^d$, there is a bi-sequence $(x_k)$ in
$\R ^d$ such that $(\rho _k *\delta _{x_k})$ is a solution of
equation (1) with noise $(\mu _k*\delta _{v_k})$ where $x
_k=v_k+\phi (x_{k-1})$.

Further, if $(\rho _k)$ is a fundamental solution of equation of (1)
with noise $(\mu _k)$, then $(\rho _k *\delta _{x_k})$ is a
fundamental solution of equation (1) with noise $(\mu _k*\delta
_{v_k})$.\el

\bo Choose $(x_k)$ in $\R ^d$ so that $x _k=v_k+\phi (x_{k-1})$
using Lemma \ref{sdc}.  Define $\lam _k = \rho _k*\delta _{x_k}$ for
any $k \in \Z$.  Then $\mu _k * \delta _{v_k}* \phi (\lam _{k-1}) =
\mu_k *\delta _{v_k} * \phi (\rho _{k-1}) *\phi (\delta _{x_{k-1}})
= \rho _{k} *\delta _{x_k} =  \lam _k$ for any $k \in \Z$.

Further, if $(\rho _k)$ is a fundamental solution of equation of (1)
with noise $(\mu _k)$ and $(\lam _k ')$ is a solution of equation
(1) with noise $(\mu _k *\delta _{v_k})$. Then by first part, $(\rho
_k')$ is a solution of equation (1) with noise $(\mu _k)$ where
$\rho _k'= \lam _k '*\delta _{-x_k}$.  Since $\rho _k$ is a
fundamental solution to equation (1) with noise $(\mu _k)$, there is
a $\nu \in M^1(\R ^d)$ such that $\rho _k '= \rho _k *\phi ^k(\nu
)$, hence $\lam _k '= \rho _k' *\delta _{x_k} = \rho _k
*\phi ^k(\nu )*\delta _{x_k} = \lam _k *\phi ^k(\nu )$.  Thus,
$(\lam _k)$ is a fundamental solution to equation (1) with noise
$(\mu _k *\delta _{v_k})$. \eo

We now prove the main result of this section.

\bo $\!\!\!\!\!${\bf ~of Theorem \ref{sdfs}~~~} Suppose there is a
strongly $\tau$-decomposable measure $\lam \in M^1((\R ^d)^\Z)$ with
co-factor $\mu$.  Then $\lam = \mu *\tau (\lam )$ and $\tau ^n(\lam
) \to \delta _0$.  Define $\mu ^{(n)}= \mu
*\tau (\mu )*\cdots *\tau ^{n-1}(\mu )$ for any $n \geq 1$.  Then
$\lam = \mu ^{(n)} *\tau ^n(\lam)$ for any $n \geq 1$.  Since $\tau
^n(\lam ) \to \delta _0$, we get that $\mu ^{(n)}\to \lam$.  Let
$\lam _k = P_k(\lam )$.  Then it follows easily that $(\lam _k)$ is
a solution of equation (1).  Take any other solution $(\lam _k')$ of
equation (1).  By Kolmogorov consistency theorem, there is a $\lam
'\in M^1((\R ^d)^\Z)$ such that $P_k(\lam ')= \lam _k' $.  By the
uniqueness part of Kolmogorov consistency theorem, $\lam ' = \mu
*\tau (\lam ')$: note that $P_k\tau = \phi P_{k-1}$. This implies
that $\lam '= \mu ^{(n)}*\tau ^n(\lam ')$ for any $n \geq 1$.  Since
$\mu ^{(n)} \to \lam$, $(\tau ^n(\lam '))$ is relatively compact.
Let $\rho \in M^1((\R ^d)^\Z)$ be a limit point of $(\tau ^n(\lam
'))$. Then $\tau ^k(\rho )$ is also a limit point of $(\tau ^n(\lam
'))$, hence $\lam '= \lam *\rho = \lam *\tau ^k(\rho)$ for any $k
\in \Z$. Let $\rho _0 = P_0(\rho)$.  Then for any $k \in \Z$, $\lam
'= \lam * \tau ^k(\rho )$ implies that $\lam _k'= \lam _k*\phi
^k(\rho _0)$. This proves that $(\lam _k)$ is a fundamental
solution.

Suppose the equation (1) has a solution.  By lemma \ref{sft}, there
are $w_k$ in $\R ^d$ and $\varrho _k\in M^1(\R ^d)$ such that $\prod
_{i=0}^{n-1} \phi ^i(\mu _{k-i}) * \phi ^{n}(\delta _{w_{n-k}}) \to
\varrho _k$ as $n \to \infty$ for any $k <0$.  Let $v_k = \phi
(w_{-k+1})-w_{-k}$ for $k<0$.  Then
$$\begin{array}{ll} \sum _{i=0}^{n-1} \phi ^i(v_{k-i}) & = \sum _ {i=0}^
{n-1} \phi ^i [\phi (w_{-k+i+1})-w_{-k+i}] \\ & = \sum _{i=0}^{n-1}
\phi ^{i+1}(w_{i-k+1}) - \sum _{i=0}^{n-1} \phi ^i (w_{i-k})\\ & =
\phi ^{n} (w_{n-k}) - w_{-k} \end{array}$$ for any $k <0$ and $n
\geq 1$.

Let $v_k =0$ for all $k \geq 0$ and $\mu _k '= \mu _k *\delta
_{v_k}$ for all $k \in \Z$.  Then $$\prod _{i=0}^{n-1} \phi ^i(\mu
_{k-i} ') = \prod _{i=0}^{n-1} \phi ^i(\mu _{k-i} )* \phi
^{n}(\delta _{w_{n-k}})
*\delta _{w_{-k}}$$ for all $k <0$.  Thus, $\prod _{i=0}^n \phi ^i(\mu
_{k-i} ') \to \varrho _k *\delta _{w_{-k}}$ for all $k <0$.

Let $\rho _k = \varrho _k *\delta _{w_{-k}}$ for $k <0$ and $\rho _k
= \prod _{i=0} ^k \phi ^i(\mu _{k-i})*\phi ^{k+1}(\rho _{-1})$ for
$k \geq 0$.  Then $\prod _{i=0 }^n \phi ^i (\mu _{k-i}') \to \rho
_k$ as $n \to \infty$ for all $k\in \Z$. This implies by Proposition
\ref{suff} that $\rho = (\rho _k)\in M^1((\R ^d)^\Z)$ is a strongly
$\tau$-decomposable measure with co-factor $\mu '= \mu*\delta _v$
where $v= (v_k)$. Now second part of (ii) follows from (i) and Lemma
\ref{seqn}. \eo

The following improves Theorem \ref{sdfs} for symmetric noise which
could be compared with factorization theorem of \cite{Si2}.

\bc\label{sn} Let $(\mu _k)$, $\phi$ and $\tau$ be as in Theorem
\ref{sdfs}. Suppose each $\mu _k$ is symmetric.  Then equation (1)
has a solution if and only if there is a symmetric strongly
$\tau$-decomposable measure $\lam \in M^1((\R ^d)^\Z)$ with
co-factor $\mu = (\mu _k)$. \ec

\bo It is sufficient to prove the only if part.  Suppose equation
(1) has a solution.  Then by Theorem \ref{sdfs}, there is a $v\in
(\R^d )^\Z$ and $\rho \in M^1((\R ^d)^\Z)$ such that $\rho$ is
strongly $\tau$-decomposable with co-factor $\mu
*\delta _v$ where $\mu = (\mu _k) \in M^1((\R ^d)^\Z)$.  This
implies that $\mu *\delta _v = \mu *\delta _v *\tau (\rho ) = \mu
^{(n)}*\delta _{v^{(n)}}*\tau ^n(\rho )$ where $\mu ^{(n)} = \prod
_{i=0}^{n-1} \tau ^i(\mu )$ and $v^{(n)} = \sum _{i=0}^{n-1} \tau
^i(v )$ for all $n \geq 1$.  Since each $\mu _k$ is symmetric, by
the uniqueness part of Kolmogorov consistency theorem, we get that
$\mu$ is symmetric and hence each $\mu ^{(n)}$ is symmetric.  By
Lemma \ref{sym}, we get that $(\mu ^{(n)})$ converges to $\lam \in
M^1((\R ^d)^\Z)$.  Since $\mu ^{(n+1)} = \mu * \tau (\mu ^{(n)})$,
we get that $\lam = \mu *\tau (\lam )$, hence $\lam = \mu
^{(n)}*\tau ^n(\lam )$.  Since $\mu ^{(n)}\to \lam$, $(\tau ^n(\lam
))$ is relatively compact and for any limit point $\nu$ of $(\tau
^n(\lam ))$, we have $\lam = \lam *\nu$.  This implies that $\nu =
\delta _0$, hence $\tau ^n(\lam ) \to \delta _0$.  This proves that
$\lam$ is a symmetric strongly $\tau$-decomposable measure with
co-factor $\mu$. \eo

As a consequence of Theorem \ref{sdfs} we now relax the requirement
for a solution to be fundamental.

\bp Let $(\mu _k)$ and $\phi$ be as in Theorem \ref{sdfs}.  Let
$(\lam _k)$ be a solution of equation (1).  Then $(\lam _k)$ is a
fundamental solution if and only if to each solution $(\nu _k)$
there is a $\nu \in M^1(\R ^d)$ such that $\nu _0 =\lam _0
*\nu$.\ep

\bo It is sufficient to prove the "if" part.  Since $(\lam _k)$ is a
solution of (1), Theorem \ref{sdfs} implies that equation (1) has a
fundamental solution $(\rho _k)$.  Then there is a $\lam \in M^1(\R
^d)$ such that $\lam _k = \rho _k*\phi ^k(\lam )$ for all $k \in
\Z$. Suppose to each solution $(\nu _k)$ there is a $\nu \in M^1(\R
^d)$ such that $\nu _0 =\lam _0 * \nu$.  Then for $(\rho _k)$, there
is a $\rho \in M^1(\R ^d)$ such that $\rho _0 = \lam _0*\rho$. Thus,
$\rho _0 = \lam _0*\rho = \rho _0 *\lam *\rho$.  This implies that
$\lam *\rho = \delta _0$ and hence $\lam = \delta _v$ for some $v\in
\R ^d$. Therefore, $\lam _k = \rho _k *\phi ^k(\delta _v)$. Since
$(\rho _k)$ is a fundamental solution, $(\lam _k)$ is also a
fundamental solution. \eo

The next proof explores the uniqueness of fundamental solution.

\bo $\!\!\!\!\!${\bf ~of Proposition \ref{ufs}~~~} There exist
$\rho$ and $\lam$ in $M^1(\R ^d)$ such that $\rho _k = \lam _k*\phi
^k(\rho)$ and $\lam _k = \rho _k*\phi ^k(\lam )$ for all $k$. This
implies that $\lam _0 = \lam _0 * \rho *\lam $, hence $\rho *\lam =
\delta _0$.  Thus, $\rho = \delta _v$ for some $v\in \R ^d$.  \eo

Now the proof of the relation between fundamental solutions and
extremal solutions.

\bo $\!\!\!\!\!${\bf ~of Theorem \ref{fes}~~~} Let $(\lam _k)$ be a
fundamental solution of (1).  Suppose there solutions $(\lam
_{1,k})$ and $(\lam _{2,k})$ of (1) such that $\lam _k = a\lam
_{1,k}+b\lam _{2,k}$ for some constants $0\leq a,b\leq 1$ with
$a+b=1$.  Since $(\lam _k)$ is a fundamental solution, there are
$\rho _1, \rho _2 \in M^1(\R ^d)$ such that $\lam _{1,k}= \lam
_k*\phi ^k(\rho _1)$ and  $\lam _{2,k}= \lam _k*\phi ^k(\rho _2)$.
This implies that $\lam _k = \lam _k*[a\phi ^k(\rho _1)+b\phi
^k(\rho _2)]$ and hence $a\phi ^k(\rho _1)+b\phi ^k(\rho _2) =\delta
_0$. Thus, $\rho _1 = \delta _0$ or $\rho _2 = \delta _0$, hence
$\lam _k = \lam _{1,k}$ or $\lam _k= \lam _{2,k}$.  This proves that
any fundamental solution is a extremal solution of (1).

Suppose $(\lam _k')$ is a extremal solution.  By Theorem \ref{sdfs},
equation (1) has fundamental solutions.  Let $(\lam _k)$ be a
fundamental solution of equation (1).  Then there is a $\rho \in
M^1(\R ^d)$ such that $\lam _{k} '= \lam _k
*\phi ^k(\rho )$.  Let $x\in \R ^d$ be a point in the support of
$\rho$. Fix $n \in \N$, let $B_n = \{ v\in \R ^d\ \mid
||v-x||<{1\over n} \}$ and $a_n = \rho (B_n)$. Then $a_n
>0$. Define $\sigma _n = {1_{B_n}\rho \over a_n}$ and if $b_n =
1-a_n
>0$, define $\sigma _n' = {\rho - a_n\sigma _n \over b_n}$ otherwise
define $\sigma _n'=\delta _0$.  Then $\sigma _n$ and $\sigma _n'$
are in $M^1(\R^d)$ and $\rho = a_n \sigma _n +b_n \sigma _n'$. This
implies that $\lam _k '= a_n (\lam _k*\phi ^k(\sigma _n) )+b_n (\lam
_k*\phi ^k(\sigma _n'))$.  Since $\lam _k$ is a solution, $(\lam
_k*\phi ^k(\sigma _n))$ and $(\lam _k*\phi ^k(\sigma _n'))$ are also
solutions.  Thus, $(\lam _k')$ is written as a convex combination of
two solutions.  Since $(\lam _k')$ is a extremal solution, we get
that $\lam _k'=\lam _k *\phi ^k(\sigma _n)$ or $\lam _k'= \lam
_k*\phi ^k(\sigma _n ')$.  If $\lam _k'= \lam _k*\phi ^k(\sigma
_n')$, then $\lam _k'= a_n (\lam _k*\phi ^k(\sigma _n))+b_n \lam
_k'$, hence $\lam _k '=\lam _k*\phi ^k(\sigma _n)$. Thus, in any
case we have $ \lam _k '=\lam _k*\phi ^k(\sigma _n)$ for all $k \in
\Z$ and all $n \geq 1$.  Since $\sigma _n = {1_{B_n}\rho \over
a_n}$, $\sigma _n \to \delta _x$.  Thus, $\lam _k'= \lam _k*\phi
^k(\delta _x)$.  This proves that any extremal solution is also a
fundamental solution.

Let $\cE$ be the set of all extremal solutions of equation (1).  If
$\cE \not = \emptyset$, then by the second part we get that all
solutions in $\cE$ are fundamental solutions. Fix a solution $(\lam
_k)$ in $\cE$.  For any $v\in \R ^d$, define $\lam _{k,v}= \lam
_k*\phi ^k(\delta _v)$.  Then $(\lam _{k,v})$ is also a fundamental
solution, hence by the first part $(\lam _{k,v})$ is in $\cE$.  By
Proposition \ref{ufs}, the map $x\mapsto (\lam _{k,v})$ is onto.
Suppose for some $v,w$, $(\lam _{k,v})= (\lam _{k,w})$.  Then $\lam
_0 *\delta _v = \lam _0*\delta _w$, hence $\lam _0 *\delta _{v-w}=
\lam _0$. This implies that $v=w$.  Thus the correspondence
$v\mapsto (\lam _{k,v})$ between $\R ^d$ and $\cE$ is bijective. \eo

Using results proved here we show that fundamental solutions are
also infinitely divisible provided the noise process $(\mu _k)$
consists of infinitely divisible measures.

\bp\label{infd}  Let $(\mu _k)$ and $\phi$ be as in Theorem
\ref{sdfs} and $(\lam _k)$ be a fundamental solution of (1). Suppose
each $\mu _k$ is infinitely divisible.  Then each $\lam _k$ is also
infinitely divisible.\ep

\bo By Kolmogorov consistency theorem there is a $\mu \in M^1((\R
^d)^\Z)$ such that $P_k(\mu ) = \mu _k$.  Let $(\lam _k)$ be a
fundamental solution of equation (1) and $\tau$ be defined as in
equation (L).  Then by Theorem \ref{sdfs}, there is a $v\in (\R
^d)^\Z$ and $\rho \in M^1((\R ^d)^\Z)$ such that $\rho$ is strongly
$\tau$-decomposable with co-factor $\mu *\delta _v$.  This implies
that $\mu
* \tau (\mu )*\cdots *\tau ^n(\mu ) *\delta _{v^{(n)}} \to \rho $
where $v^{(n)}= \sum _{k=0}^n \tau ^k(v)$.  For any $k \in \Z$,
considering the $k$-th projection, we get that $\mu _k *\phi (\mu
_{k-1})*\cdots *\phi ^n(\mu _{k-n})
*\delta _{x_k} \to \rho _k$ where $x_k = P_k (v^{(n)})$ and
$\rho _k = P_k(\rho )$.   Suppose each $\mu _k$ is infinitely
divisible.  Since the set of infinitely divisible measures in
$M^1(\R ^d)$ is a closed set, closed under convolution product, we
get that $\rho _k$ is infinitely divisible.  Since $\rho $ is
strongly $\tau$-decomposable with co-factor $\mu *\delta _v$ and
$P_k(\rho ) = \rho _k$, we get that $(\rho _k)$ is a fundamental
solution of equation (1) for the noise $(\mu _k*\delta _{v_k})$
where $v_k = P_k(v)$.  By Lemma \ref{seqn}, $(\rho _k *\delta
_{w_k})$ is a fundamental solution of equation (1) for the noise
$(\mu _k)$.  Since $(\lam _k )$ is a fundamental solution of
equation (1) for the noise $(\mu _k)$, by Proposition \ref{ufs},
there are $y_k \in \R ^d$ such that $\lam _k=\rho _k
*\delta _{w_k}*\delta _{y_k}$.  Since each $\rho _k$ is infinitely
divisible, each $\lam _k=\rho _k *\delta _{w_k}*\delta _{y_k}$ is
also infinitely divisible.  \eo

\end{section}

\begin{section}{Noise is stationary}

We now look at the situation when the noise process $(\xi _k)$ is
stationary with common law $\mu$.   Thus, equation (1) becomes
$$\lam _k = \mu *\phi (\lam _{k-1}), ~~~~~ k \in \Z \eqno (4)$$ where
$\phi$ is a linear map on $\R ^d$.  In this case the well known
contraction subspace associated to a linear map $\phi$ plays a
crucial role: recall that contraction subspace of a linear map
$\phi$ is denoted by $C(\phi )$ and is defined by
$$C(\phi ) = \{ v\in \R ^d \mid \phi ^n(v) \to 0 ~~{\rm as}~~ n \to
\infty \}$$ and it is easy to see that $C(\phi )$ is a subspace of
$\R ^d$.

In case noise is stationary we obtain the following.

\bt\label{sn} Suppose the noise process $(\xi _k)$ is stationary
with common law $\mu$ and $\phi$ is an invertible linear map on $\R
^d$.  Then the following are equivalent:

\be \item equation (4) has a (fundamental) solution;

\item the noise law $\mu$ is supported on a coset $u+C(\phi )$ of
$C(\phi )$ for some $u \in \R ^d$ and satisfies the following
logarithmic moment condition
$$\int \log (||v||+1)d \mu (v) <\infty .$$ \ee In this case,
$(\nu *\delta _{u_k})$ is a fundamental solution where $(u_k)$ is
given by Lemma \ref{seqn} so that $u_k = u+\phi (u_{k-1})$ and $$\nu
= \mu *\delta _{-u}*\phi (\nu ) = \lim _{k\to \infty} \prod _{i=1}^k
\phi ^{i-1}(\mu ) *\delta _{-u_k}$$ in other words, $\nu$ is a
strongly $\phi$-decomposable measure with co-factor $\mu*\delta
_{-u}$. \et

We first establish the relevance of $C(\phi )$.

\bl\label{nc} Suppose $C(\phi ) = \{ 0 \}$.  Then equation (4) has a
solution $(\lam _k)$ if and only if $\mu$ is a dirac measure. \el

\bo Suppose $(\lam _k)$ is a solution of (4).  Replacing $\lam _k$
and $\mu$ by $\lam _k*\check \lam _k$ and $\mu*\check \mu$
respectively, we may assume that $\mu$ is symmetric.  Now by
iterating the equation (4), we get for fixed $k \in \Z$ that $\lam
_k = \prod _{i=1}^{n} \phi ^{i-1}(\mu )
* \phi ^n(\lam _{k-n})$ for all $n\geq 1$.

Let $\mu _n =\prod _{i=1}^{n} \phi ^{i-1}(\mu )$.  Then $\mu _n =
\mu _k * \prod _{i=k+1}^ {n} \phi ^{i-1}(\mu)$ for all $n\geq k$,
hence it follows from Lemma \ref{sym} that $\mu _n \to \nu \in
M^1(\R ^d)$. This implies that $\nu = \mu*\phi (\nu )$, hence $\nu =
\mu _n*\phi ^n(\nu )$ for all $n \geq 1$.  Thus, $(\phi ^n(\nu ) )$
is relatively compact and any limit point $\sigma$ of $(\phi ^n(\nu)
)$ satisfies $\nu = \nu *\sigma$, hence $\sigma = \delta _0$. This
implies that $\phi ^n(\nu ) \to \delta _0$ and hence $\nu (C(\phi )
)= 1$. Since $C(\phi ) = \{0 \}$, $\nu = \delta _0$. Thus, $\mu =
\delta _0$. \eo

\bo $\!\!\!\!\!${\bf ~of Theorem \ref{sn}~~~} Suppose equation (4)
has a solution $(\lam _k)$.  Then applying Lemma \ref{nc} to
$V/C(\phi )$, we get that $\mu$ is supported on a coset $C(\phi )$.
Let $\mu _s= \mu *\check \mu$ and $\lam _{k,s} = \lam _k*\check \lam
_k$. Then $\lam _{k,s} = \mu _s*\phi (\lam _{k-1,s}) = \prod
_{i=0}^n \phi ^i (\mu _s) * \phi ^{n+1}(\lam _{k-n-1,s})$. By Lemma
\ref{sym} we get that $\prod _{i=0}^ {n} \phi ^{i}(\mu _s) \to \nu
\in M^1(\R ^d)$ as $k \to \infty$. Since $\mu$ is supported on a
coset $C(\phi )$, $\mu _s$ is supported on $C(\phi )$.  This implies
by \cite{Za} that $\int \int \log (||v-w||+1) d\mu (v) d\mu (w)
<\infty$.  By Fubini's Theorem $\int \log (||v-w||+1) d\mu (v)
<\infty$ for some $w\in \R ^d$. Since $v\mapsto \log (||v||+1)$ is
subadditive, that is $\log (1+||v_1+v_2||) \leq \log (1+||v_1||)
+\log (1+||v_2||)$, we get that $\int \log (||v||+1) d\mu (v)
<\infty$.

Assume that $\mu$ is supported on a coset $u +C(\phi )$ and $\int
\log (||v||+1) d\mu (v) <\infty$.  Then as above using subadditivity
of the map $v\mapsto \log (||v||+1)$ we get that $\int \log
(||v-u||+1) d\mu (v) <\infty$.  By \cite{Za} we get that $\prod
_{i=1}^k \phi ^{i-1} (\mu *\delta _{-u})\to \nu \in M^1(C(\phi ))$.
Then $\nu = \mu *\delta _{-u}* \phi (\nu )$.  Take $\nu _k = \nu$
and $\mu _k = \mu *\delta _{-u}$ for all $k \in \Z$.  Then $\phi
^n(\nu _{k-n}) = \phi ^n(\nu ) \to \delta _0$.  This implies that
$(\nu _k) \in M^1((\R^d)^\Z)$ is strongly $\tau$-decomposable with
co-factor $(\mu _k)\in M^1((\R^d)^\Z)$.  By Theorem \ref{sdfs},
$(\nu _k)$ is a fundamental solution to equation (1) for the noise
$(\mu _k)$. Using Lemma \ref{sdc} we may find a sequence $(u_k)$
such that $u_k = u+\phi (u_{k-1})$ for all $k \in \Z$.  Define $\lam
_k = \nu * \delta _{u_k}$ for all $k \in \Z$. Then By Lemma
\ref{seqn}, we get that $(\lam _k)$ is a fundamental solution to
equation (1) for the noise $(\mu _k*\delta _u)$, that is equation
(4) has a (fundamental) solution. \eo

\end{section}

\begin{section}{Noise with independent $\ell _p$-paths}

We now look at the situation where the noise process has independent
$\ell _p$-paths.

Recall that $\ell _p =\{ (x_k) \in (\R ^d)^\Z \mid \sum ||x_k||^p
<\infty \}$ for $1\leq p<\infty$.  Then $\ell _p$ are proper dense
subspaces of $(\R ^d)^\Z$.  But the spaces $\ell _p$ with norm
$||\cdot ||_p$ given by $||(x_k)||_p ^p = \sum ||x_k||^p$ are
separable Banach spaces for $1\leq p<\infty$.

Fix $p \in [1, \infty)$ and $\phi$ be a linear transformation on $\R
^d$.  Let $\iota \colon \ell _p \to (\R ^d)^\Z$ be the natural
inclusion and $\ap \colon \ell _p \to \ell _p$ be defined by $\ap
(v_k) = (\phi (v_{k-1}))$. It is easy to check that $\iota$ and
$\ap$ are continuous linear maps, $\tau \odot \iota = \iota \odot
\ap$.  If the noise process $(\xi _k)$ has independent $\ell
_p$-paths, then for $\mu = (\mu _k)\in M^1((\R ^d)^\Z)$, we have
$\mu (\iota (\ell _p))=1$.  So, $\mu$ can be realized as a Borel
probability measure on $\ell _p$.  We follow this notations and
realization in this section.

If the noise process $(\mu _k)$ has independent $\ell _p$-paths,
then we introduce logarithmic moment $M_p$ by

$$M_p =  \lim _{n\to \infty}
\int _{(\R ^d)^{2n+1}} \log ( (\sum _{-n\leq k \leq n}
||v_k||^p)^{1\over p}+1) d\mu _n'(v) $$ where $\mu _n'$ is the
product measure $\mu _{-n} \times \cdots \times \mu _n$ defined on
$(\R ^d)^{2n+1}$.   We now obtain the following result which is
similar to the stationary noise situation: the main ingredients in
the proof is the realization of the noise process in the Banach
space $\ell _p$ and results of \cite{Si1} and \cite{Si2}.

\bt\label{bp} Suppose the noise process $(\mu _k)$ has independent
$\ell _p$-paths and $\phi$ is an invertible linear map on $\R ^d$.
Then the following are equivalent: \be

\item there exists $(x_k)$ in $\ell _p$ such that
each $\mu _k$ is supported on the coset $C(\phi )-x_k$ and satisfies
the logarithmic moment condition that $M_p <\infty$;

\item there is a (fundamental) solution $(\lam _k)$ of equation (1) such that
$(\lam _k*\delta _{y_k})$ has independent $\ell _p$-paths for some
$(y_k)$ such that $y_k= -x_k+\phi (y_{k-1})$.

\ee

\et

\br Theorem \ref{bp} has a shift $(y_k)$ because the equation
$y_k=x_k+\phi (y_{k-1})$ may not have a solution $(y_k)$ in $\ell
_p$ for a given $(x_k)$ in $\ell _p$.  For instance, fix $b>1$ and
define $\phi \colon \R \to \R$ by $\phi (t)= b^2t$ for all $t \in
\R$.  Take $x_k = b^{-|k|}$ for all $k \in \Z$. Then for $1\leq
p<\infty$, $(x_k)\in \ell _p$ but there is no $(y_k) \in \ell _p$
such that $y_k = x_k+ \phi (y_{k-1})$ for all $k \in \Z$. Because if
there is a $(y_k) \in \ell _p$ such that $y_k = x_k+ \phi
(y_{k-1})$, then $y_0 = \sum _{i=0} ^k \phi ^i(x_{-i}) +y_{-k-1}=
\sum _{i=0}^k b^i + y_{-k-1}$ for all $k \geq 1$.  Since $(y_k)\in
\ell _p$, $y_{-k}\to 0$ as $k \to \infty$.  Thus, $y_0 = \sum _{i
\geq 0} b ^i$ which is a contradiction to $b>1$.  In general, the
equation $y=x+\tau (y)$ has a solution in $\ell _p$ for any given
$x\in \ell _p$ if and only if $1$ is not in the spectrum of $\tau$
as $\tau$ has no fixed points in $\ell _p$. \er

\bex

We now give examples of noise processes that have independent $\ell
_p$-paths and also satisfy the logarithmic moment condition.

\be

\item  Suppose $(\mu _k)$ is a bi-sequence in $M^1(\R ^d)$ such that
the corresponding absolute moments $m_k = \int ||x|| d\mu _k(x)$ are
summable, that is $\sum m_k <\infty$.  Then if $(X_k)$ is a
independent bi-sequence of random variables such that law of $X_k$
is $\mu _k$, then $\sum \int ||X_k|| <\infty$. This implies that
$\int \sum ||X_k|| <\infty$ and hence $\sum ||X_n || <\infty$ a.s.
Thus, $(\mu _k)$ has independent $\ell _1$-paths.  Since $\log
(|t|+1)\leq |t|$, we get that $\int \log (||v||+1) d\mu _k \leq
m_k$.  Using the subadditivity of $t\mapsto \log (|t|+1)$, we can
conclude that $(\mu _k)$ satisfies the moment condition of Theorem
$\ref{bp}$.

\item Fix $0<a <1$.  For $n\geq 1$, let $\omega _n$ (resp. $\omega
_{-n}$) be the uniform measure defined on the closed interval
$[0,a^n]$ (resp. on $[-a^n, 0]$) and $\omega _n '$ (resp. $\omega
_{-n}'$) be the uniform measure defined on the closed interval
$[n,n+1]$ (resp. on $[-n-1, -n]$). Suppose $(X_n)$ is a bi-sequence
of independent random variables such that the law of $X_n$ is given
by
$$\mu _n = \left \{
\begin{array}{ll}
(1-a^{|n|}) \omega _n +a^{|n|} \omega _n'& n\not =0 \\
\delta _0 & n=0 \\
\end{array}\right.  $$

for all $n\in \Z$.  Then $$\sum P(X_n^2 >1) = \sum \mu _n(\R
\setminus [-1,1]) = \sum _{n\not =0} a^{|n|} <\infty ,$$
$$\sum E(X_n ^21_{\{|X_n|^2 <1\}}) = \sum (1-a^{|n|})\int _{-1} ^1 x^2
d\omega _n  (x) \leq \sum a^{2|n|} <\infty$$ and
$$\sum |\Var (X_n^21_{\{|X_n|^2 <1\}})| \leq 2\sum a^{4|n|}
<\infty .$$ By Kolmogorov three series theorem (cf. Theorem 5.3.3 of
\cite{Ch}), we get that $\sum X_n^2$ converges a.s.  So the laws
$(\mu _k)$ have independent $\ell _2$-paths.  This noise process
also verifies the logarithmic moment condition given in Theorem
\ref{bp}. Similar examples can be obtained using Kolmogorov three
series theorem.

\ee

\eex

We first compare $||\ap ||_p$ and $||\phi ||$.

\bl\label{nr} $||\ap ||_p= ||\phi ||$. \el

\bo For $v\in \ell _p$, we have $||\ap (v)||_p \leq ||\phi ||
||v||$.  This implies that $||\ap ||_p \leq ||\phi ||$.  For any
$v_0 \in \R ^d$ such that $||v_0|| =1$, define $v_k = 0$ for all $k
\not =0$ and take $v= (v_k)$.  Then $v\in \ell _p$ with $||v||_p =1$
and $||\ap (v)|| _p= ||\phi (v_0)||$, hence $||\phi ||\leq ||\ap
||_p$. This proves that $||\ap ||_p= ||\phi ||$. \eo

We next show the relevance of $C(\phi )$.

\bl\label{bnr1} Assume that the noise process $(\mu _k)$ has
independent $\ell _p$-paths.  If $(\lam _k)$ is a solution of (1)
that has independent $\ell _p$-paths, then there exists a $(x_k) \in
\ell _p$ such that each $\mu _k$ is supported on the coset
$x_k+C(\phi )$. \el

\bo We first show that $C(\ap ) \subset \prod C(\phi )$. For
$v=(v_k) \in \ell _p$, $$||\ap ^n(v)|| = (\sum _k ||\phi
^n(v_{k-n})||^p )^{1\over p}= (\sum _k ||\phi ^n(v_{k})||^p
)^{1\over p}$$ for all $n$. So, if $\ap ^n (v) \to 0$, then $||\phi
^n(v_k)||\to 0 $ for all $k \in \Z$, hence $v_k \in C(\phi )$.

Suppose equation (1) has a solution $(\lam _k)$ that has independent
$\ell _p$-paths.  Then by Kolmogorov consistency theorem and by its
uniqueness part we get that there are $\mu$ and $\lam $ in $M^1(\ell
_p)$ such that $P_k(\mu )=\mu _k$, $P_k(\lam )=\lam _k$ and $\lam =
\mu *\ap (\lam )$.  Replacing $\mu $ by $\mu
*\check \mu$, we may assume that $\mu$ is symmetric and $\mu$ is a
co-factor of a $\ap$-decomposable symmetric measure. By 3.3 of
\cite{Si2}, there is a strongly $\ap$-decomposable measure $\nu$ on
$\ell _p$ with co-factor $\mu$.  Now the result follows from
Corollary 3.3 of \cite{Si1}. \eo

\bo $\!\!\!\!\!${\bf ~of Theorem \ref{bp}~~~}  Assume $(\mu _k)$ has
independent $\ell _p$-paths.  By Kolmogorov consistency theorem
there is a $\mu \in M^1(\ell _p)$ such that $P_k(\mu ) = \mu _k$. We
first observe the following:

$$\begin{array}{ll}\int \log (||v||_p + 1) d\mu (v) & = \int \log
((\sum _k ||v_k|| ^p)^{1\over p}+1) d\mu (v) \\ & = \lim _{n \to
\infty}\int _{(\R ^d)^{2n+1}}   \log ((\sum _{-n\leq k \leq n}
||v_k||^p)^{1\over p} +1) d\mu _n' \end{array} \eqno {\rm (M_p')}$$
for $p\in [1, \infty )$  where $\mu _n'$ is the product measure $\mu
_{-n} \times \cdots \times \mu _n$ defined on $(\R ^d)^{2n+1}$.

Let $X_0 = \{ (v_k) \in \ell _p \mid v_k \in C(\phi ) \}$.  Then
$X_0$ is a closed subspace of $\ell _p$ and $\ap (X_0) =X_0$. Since
$||\ap ||_p= ||\phi ||$, it could easily be verified that the
spectral radius of $\ap$ restricted to $X_0$ is equal to the
spectral radius of $\phi$ restricted to $C(\phi )$ which is less
than one unless $C(\phi )= \{ 0 \}$.

In order to prove the result, we assume that $C(\phi )$ is
nontrivial.  Thus, spectral radius of $\ap$ restricted to $X_0$ is
less than one.

Suppose there is a $x=(x_k)\in \ell _p$ such that $\mu _k$ is
supported on $C(\phi )-x_k$ and has the logarithmic moment condition
that $M_p <\infty$. Then using $(M_p')$ and the subadditivity of
$v\mapsto \log (||v||_p+1)$ we get that
$$ \int \log (||v+x||_p+1) d\mu (v) \\
\leq \int \log ( ||v||_p +1) d\mu (v) + \log (||x||_p+1) < \infty.
$$  Since $\mu _k$ is supported on $C(\phi )-x_k$, $\mu * \delta _x
(X_0)=1$.  Since $\ap$ restricted to $X_0$ has spectral radius less
than one, 1.5 of \cite{Si1} implies that $\prod _{i=0}^{n-1}\phi
^i(\mu *\delta _x) \to \lam \in M^1(X_0)$.  This implies that $\lam
= \mu
*\delta _x *\ap (\lam )$ and $\ap ^n(\lam )\to \delta _0$.  It follows
from Theorem \ref{sdfs} that $P_k(\lam )$ is a fundamental solution
to equation (1) for the noise $(\mu _k *\delta _{x_k})$. By Lemma
\ref{seqn}, equation (1) has a fundamental solution for the noise
$(\mu _k)$.

Suppose equation (1) has a solution $(\lam _k)$ such that $(\lam
_k*\delta _{z_k})$ has independent $\ell _p$-paths.  Then $\lam _k =
\mu _k *\phi (\lam _{k-1})$ implies $$\lam _k*\delta _{z_k} = \mu _k
*\phi (\lam _{k-1})*\delta _{z_k} = \mu _k*\delta _{a_k} *\phi
(\lam _{k-1}*\delta _{z_{k-1}})$$ where $a_k = z_k-\phi (z_{k-1})$.
Since $\mu _k$, $\lam _k*\delta _{z_k}$ have independent $\ell
_p$-paths, $(a_k)$ is in $\ell _p$.  Since $(\mu_k)$ has independent
$\ell _p$-paths, $(\mu _k*\delta _{a_k})$ also has independent $\ell
_p$-paths. Applying Lemma \ref{bnr1} to $\lam _k*\delta _{z_k} = \mu
_k*\delta _{a_k} *\phi (\lam _{k-1}*\delta _{z_{k-1}})$, we get that
each $\mu _k$ is supported on $x_k+C(\phi )$ for some $x =(x_k)\in
\ell _p$. By Kolmogorov consistency theorem and by its uniqueness
there is a $\lam \in M^1(\ell _p)$ such that $\lam = \mu
*\delta _a* \ap (\lam )$ where $P_k(\lam ) = \lam _k*\delta _{z_k}$
and $a=(a_k)$.   This implies that $\lam *\check \lam = \mu *\check
\mu * \ap (\lam *\check \lam)$.  By Theorem 3.3 of \cite{Si2} we get
that $\prod _{i=1}^ {n} \ap ^{i-1}(\mu *\check \mu) \to \nu \in
M^1(\ell _p)$ as $n \to \infty$.  Since each $\mu _k$ is supported
on $x_k+C(\phi )$, $\mu _k *\check \mu _k$ is supported on $C(\phi
)$.  So $\mu *\check \mu$ is supported on $X_0$.  This implies by
1.5 of \cite{Si1} that $\int \int \log (||v-w||_p+1) d\mu (v) d\mu
(w) <\infty$.  By Fubini's theorem $\int \log (||v-w||_p+1) d\mu (v)
<\infty$ for some $w\in \R ^d$. Since $v\mapsto \log (||v||_p+1)$ is
subadditive, we get that $\int \log (||v||_p+1) d\mu (v) <\infty$.
Now the result follows from $(M_p')$.

\eo

\end{section}

\begin{section}{Remarks}

We now make a few remarks.  The first one provides a counter-example
to show that having moment condition on the individual marginals of
noise process as in Theorem \ref{sn} is not sufficient for existence
of (fundamental) solutions.

(1) Let $\mu _0$ be the uniform measure supported on $[0,1]$ and
$0<a<1$.  Define $$\mu _k = \phi ^k(\mu _0 )$$ for all $k \in \Z$
where $\phi (t) = at$ for all $t\in \R$.  Then $$\int \log (|t|+1)
d\mu _k(t) =\int _0 ^1 \log (|a^kt|+1)dt \leq \log (a^k+1) <\infty$$
for all $k$.  But the equation $\lam _k = \mu _k*\phi (\lam
_{k-1})$, $k\in \Z$ has no solution.  Because, if there is a
solution $(\lam _k)$, then $\lam _0= \mu _0*\phi (\lam _{-1}) = \mu
_0 *\phi (\mu _{-1})*\phi ^2(\lam _{-2}) = \cdots = \prod _{k=0}^n
\phi ^k(\mu _{-k}) *\phi ^{n+1}(\lam _{-n-1}) = \mu ^{n+1}*\phi
^{n+1}(\lam _{-n-1})$ for all $n \geq 1$.  This implies that $(\mu
^{n+1}*\delta _{x_n})$ is relatively compact for some sequence
$(x_n)$ in $\R$, hence $\mu$ is a dirac measure which is a
contradiction.  Here, Theorem \ref{bp} can not be applied as it can
be seen by using Kolmogorov three series theorem that $(\mu_k)$ does
not have independent $\ell_p$-paths for any $p \in [1,\infty)$.

(2) {\it Gaussian noise:}  Let $(\mu _k)$ be a bi-sequence in
$M^1(\R ^d)$ and $\phi$ be a linear transformation on $\R ^d$.
Suppose each $\mu _k$ is Gaussian with covariance operator $A_k$. We
now claim that equation (1) has a solution if and only if $\sum
_{i=0} ^\infty (\phi ^i) ^* A_{-i} \phi ^i$ converges where $(\phi ^
i)^*$ denotes the adjoint of $\phi ^i$.  Applying Lemma \ref{seqn},
we may assume that $\mu _k$ is symmetric. Suppose equation (1) has a
solution. Define $\tau$ by equation (L). By Corollary \ref{sn},
there is a symmetric strongly $\tau$-decomposable measure $\lam \in
M^1((\R ^d)^\Z)$ with co-factor $\mu = (\mu _k)$. Then $\lam = \lim
_n (\prod _{k=0}^{n-1}\tau ^k( \mu ))$. Let $\lam _k = P_k( \lam )$.
Then by Theorem \ref{sdfs}, we get that $(\lam _k)$ is a fundamental
solution of equation (1) and $\lam _k = \lim _n \prod _{i=0}^{n-1}
\phi ^i(\mu _{k-i})$.  By considering the corresponding
characteristic functions, we get that $\sum _{i=0} ^\infty (\phi ^i)
^* A_{k-i} \phi ^i$ converges and $\lam _k$ is Gaussian with
covariance operator $\sum _{i=0} ^\infty (\phi ^i) ^* A_{k-i} \phi
^i$.  Conversely, if $\sum _{i=0} ^\infty (\phi ^i) ^* A_{-i} \phi
^i$ converges, then $\sum _{i=0} ^\infty (\phi ^i) ^* A_{k-i} \phi
^i = (\phi ^k)^* (\sum _{i=-k} ^\infty (\phi ^i) ^* A_{-i} \phi
^i)\phi ^k$ converges for all $k \in \Z$.  By taking $\lam _k$ to be
the Gaussian with covariance operator $B_k= \sum _{i=0} ^\infty
(\phi ^i) ^* A_{k-i} \phi ^i$, it may be easily verified that $\lam
_k = \lim \prod _{i=0}^{n-1} \phi ^i (\mu _{k-i})$.  Thus, by
Proposition \ref{suff} and Theorem \ref{sdfs} we get that $(\lam
_k)$ is a fundamental solution.  A similar result may be obtained
for any noise consisting of infinitely divisible distribution using
Levy-Khinchin representation (cf. 5.7 of \cite{Li}) but for
simplicity we considered the Gaussian case.

(3) {\it Gaussian and Poisson solution:}  Suppose the equation $\lam
_k =\mu _k*\phi (\lam _{k-1})$ has a solution $(\lam _k)$ consisting
of Gaussian measures $\lam _k$.  Then it follows from Cram\'er's
theorem that each $\mu _k$ is Gaussian (see \cite{LO}).  Thus,
equation (1) has a solution consisting of Gaussian measures only if
each $\mu _k$ is Gaussian. Similarly, we may conclude that equation
(1) has a solution consisting of Poisson measures only if each $\mu
_k$ is Poisson-use Raikov's theorem (see \cite{LO}).  Using similar
idea it can also be seen that if there is a solution of (1)
consisting of Gaussian (resp. Poisson) measures, then there is a
fundamental solution of (1) consisting of Gaussian (resp. Poisson)
measures.

(4) {\it Gaussian noise with no solution:} If we take $\mu _0$ in
(1) to be the Gaussian $N(0,1)$ and define $\phi$ and $\mu _k$ as in
(1).  Then $\mu _k$ is the Gaussian measure $N(0,a^{2k})$.  As in
(1), it may be verified that the equation  $\lam _k = \mu _k*\phi
(\lam _{k-1})$, $k\in \Z$ has no solution.  It may also easily be
seen that $\mu _k$ does not satisfy the condition in (2).

(5) {\it General noise:}  We now provide an example to show that in
general one may get a (fundamental) solution of equation (1) even if
$C(\phi )= \{ 0 \}$.  Let $\phi$ be a linear transformation on $\R
^d$.  Take $\mu _k$ to be the Gaussian with covariance operator
$A_k$ with $||A_k||= a^k||\phi ||^{2k}$ for some $a>1$.  Then
$$||\sum _{i=0} ^{n} (\phi ^i)^* A_{-i}\phi ^i|| \leq \sum _{i=0}^n
a^{-i} ||\phi ||^{2i} ||\phi ||^{-2i} = \sum _{i=0} ^{n}
a^{-i}<\infty $$ as $a>1$.  Thus, by (2) equation (1) has a solution
for the noise $(\mu _k)$.

\end{section}

\bigskip\medskip
\advance\baselineskip by 2pt
\begin{tabular}{ll}
C.\ R.\ E.\ Raja \\
Stat-Math Unit \\
Indian Statistical Institute (ISI) \\
8th Mile Mysore Road \\
Bangalore 560 059, India.\\
creraja@isibang.ac.in
\end{tabular}


\begin{thebibliography}{JRT-94}

\footnotesize


\bibitem [AkUY-08]{AUY} J. Akahori, C. Uenishi and K. Yano, Stochastic
equations on compact groups in discrete negative time. Probab.
Theory Related Fields 140 (2008), 569--593.

\bibitem [Ch-01]{Ch} K. L. Chung, A course in probability theory, Third
edition. Academic Press, Inc., San Diego, CA, 2001.

\bibitem [Cs-66]{Cs} I. Csisz\'ar, On infinite products of random elements
and infinite convolutions of probability distributions on locally
compact groups, Z. Wahrscheinlichkeitstheorie und Verw. Gebiete 5
(1966), 279--295.

\bibitem [Li-86]{Li} W. Linde, Probability in Banach spaces—stable and infinitely
divisible distributions. Second edition. A Wiley-Interscience
Publication. John Wiley \& Sons, Ltd., Chichester, 1986.

\bibitem [LiO-77]{LO} Ju. V. Linnik and I. V. Ostrovs'kii,
Decomposition of random variables and vectors. Translated from the Russian.
Translations of Mathematical Monographs, Vol. 48.
American Mathematical Society, Providence, R. I., 1977.

\bibitem [HiY-10]{HiY} T. Hirayama and K. Yano, Extremal solutions for
stochastic equations indexed by negative integers and taking values
in compact groups, Stochastic Process. Appl. 120 (2010), 1404--1423.

\bibitem [Ke73]{Ke} H. Kesten, Random difference equations and renewal
theory for products of random matrices, Acta Math. 131 (1973),
207--248.

\bibitem [Pa-67]{Pa} K. R. Parthasarathy, Probability measures on
metric spaces, Probability and Mathematical Statistics, No. 3
Academic Press, Inc., New York-London 1967.

\bibitem [Ra-12]{Ra} C. R. E. Raja, A stochastic difference equation
with stationary noise on groups, Canadian Journal of Mathematics 64
(2012), 1075--1089

\bibitem [Si-91]{Si1} E. Siebert, Strongly operator-decomposable
probability measures on separable Banach spaces, Math. Nachr. 154
(1991), 315--326.

\bibitem [Si-92]{Si2} E. Siebert, Operator-decomposability of Gaussian
measures on separable Banach spaces, J Theor. Prob. 5 (1992),
333--347.

\bibitem [Ta-09]{Ta} Y. Takahashi, Time evolution with and without remote
past, Advances in discrete dynamical systems, 347--361, Adv. Stud.
Pure Math., 53, Math. Soc. Japan, Tokyo, 2009.

\bibitem [To-69]{To} A. Tortrat, Convolutions d\'enombrables \'equitendues
dans un groupe topologique X, Les probabilit\'es sur les structures
alg\'ebriques (Actes Colloq. Internat. CNRS, No. 186,
Clermont-Ferrand, 1969), 327–-343. \'Editions Centre Nat. Recherche
Sci., Paris, 1970.

\bibitem [Ts-75]{Ts} B. S. Cirel'son, An example of a stochastic differential
equation that has no strong solution. (Russian) Teor. Verojatnost. i
Primenen. 20 (1975), 427--430.

\bibitem [Yo-92]{Yo} M. Yor, Tsirel'son's equation in discrete time,
Probab. Theory Related Fields 91 (1992), 135--152.

\bibitem [Za-77]{Za} O. Zakusilo, Some properties of random vectors of the
form $\sum _0 ^\infty Aî \xi _i$, Theory Probab. Math. Statist. 13
(1977), 62–-64.



\end{thebibliography}
\end{document}